\documentclass{amsart}
\usepackage{graphicx}
\usepackage{amscd}
\usepackage{amsmath}
\usepackage{amsthm}
\usepackage{amsfonts}
\usepackage{amssymb}
\usepackage{mathrsfs}
\usepackage{enumerate}
\usepackage[ pdfstartview=FitB]{hyperref}

\newcommand{\D}{\operatorname{\mathbb{D}}}
\newcommand{\N}{\operatorname{\mathbb{N}}}

\newcommand{\C}{\operatorname{\mathbb{C}}}

\newcommand{\B}{\operatorname{\mathcal{B}}}

\newcommand{\e}{\operatorname{\varepsilon}}

\newcommand{\ol}{\overline }
\DeclareMathOperator{\id}{Id}

\newtheorem{lemma}{Lemma}[section]
\newtheorem{theorem}[lemma]{Theorem}
\newtheorem{proposition}[lemma]{Proposition}
\newtheorem{corollary}[lemma]{Corollary}
\theoremstyle{definition}

\begin{document}
\author{Rapha\"el Clou\^atre}
\address{Department of Mathematics, Indiana University, 831 East 3rd Street,
Bloomington, IN 47405}
\email{rclouatr@indiana.edu}
\title[Similarity results for operators of class $C_0$]{Similarity results for operators of class $C_0$}
\subjclass{Primary 47A45; Secondary 30E05.}
\begin{abstract}
If $T$ is a multiplicity-free contraction of class $C_0$ with
minimal function $m_T$, then it is quasisimilar to the Jordan block
$S(m_T)$. In case $m_T$ is a Blaschke product with simple roots forming a Carleson sequence, we show that the relation between $T$ and $S(m_T)$ can
be strengthened to similarity. Under the additional assumption that
$u(T)$ has closed range for every inner divisor $u\in H^\infty$ of
$m_T$, the result also holds in the more general setting where the roots have bounded multiplicities.
\end{abstract}
\thanks{Research supported by NSERC (Canada) }

\maketitle

\section{Introduction}
Let $T$ be a bounded linear operator on some Hilbert space. Assume that $T$ is of class $C_0$ and admits a cyclic vector. In this case, we know that $T$ is quasisimilar to the Jordan block $S(\theta)$ where $\theta \in H^\infty$ is the minimal inner function of $T$ (see Section 2 for precise definitions). The problem we consider in this paper is as follows: under the
additional assumption that $u(T)$ has closed range for every inner divisor $u\in
H^\infty$ of $\theta$, can the relation of quasisimilarity
between $T$ and $S(\theta)$ be strenghtened to \textit{similarity}?
We will motivate this question in the next section, but let us first
make a few remarks. First, it is a classical fact that the converse
of this statement holds, in fact it is known that $u(S(\theta))$ is a
partial isometry for every inner divisor $u$ of $\theta$; see
\cite{bercOTA}. Moreover, the special case of a diagonal operator gives an interpretation of Carleson's interpolation theorem in our context of minimal functions of operators of class $C_0$. Our main result shows that this problem has a positive answer in the case where
$$
\theta(z)=z^{m_0}\prod_{j\in \N}\left(\frac{\overline{\lambda_j}}{|\lambda_j|}\frac{\lambda_j-z}{1-\overline{\lambda_j}z}\right)^{m_j}
$$
with $\sup_{j}m_j<\infty$ and $\{\lambda_j\}_j$ satisfying the Carleson condition
\begin{equation}\label{unifsep}
\inf_{k\in \N} \prod_{j\neq k} \left|\frac{\lambda_{j}-\lambda_{k}}{1-\ol{\lambda_j}\lambda_k}\right|>0.
\end{equation}

The paper is organized as follows. Section 2 deals with preliminaries.
Section 3 treats the problem in the case where each $m_j$ is
equal to one. The solution hinges on the classical Carleson
interpolation theorem. In Section 4, we give a characterization of
Carleson sequences $\{\lambda_j\}_j\subset \D$ in terms of
operators of class $C_0$ with minimal function $$\prod_{j\in
\N}\left(\frac{\ol{\lambda_j}}{|\lambda_j|}\frac{\lambda_j-z}{1-\ol{\lambda_j}z}\right).$$
This motivates the assumption on the closure of the ranges. We then extend the
similarity result to the case where the multiplicities of the roots
are bounded. In addition to the statement and
proof of the main result, Section 5 contains another crucial
ingredient: an estimate concerning the norm of the similarity matrix
between nilpotent contractions and Jordan cells.

It is appropriate here to address an insightful remark communicated to us by the referee. Our main result Theorem \ref{similar2} is concerned with a Carleson sequence $\{\lambda_j\}_j \subset \D$ where each $\lambda_j$ can have multiplicity at most $M$. This is a degenerate case of the more general situation of a finite union of Carleson sequences, namely $\{\lambda_j\}_j=\bigcup_{k=1}^M \Lambda_k\subset \D$ where $\Lambda_k\subset \D$ is a Carleson sequence (with each of its elements having multiplicity one). It is a natural question whether Theorem \ref{similar2} holds in this context. Using results from  \cite{hart1}, \cite{hart2}, \cite{Vas2} and \cite{NikShift}, it is possible to generalize Proposition \ref{orthogonal} after some significant modification of our present argument. However, a full generalization of Theorem \ref{similar2} would require an extension of Proposition \ref{nilpotent} from the nilpotent case to that where the operator is merely algebraic. At this time, we feel this is the major obstacle to obtaining the stronger theorem. It is our hope that we may settle the question in future work.

The author wishes to thank his advisor Hari Bercovici for his generosity, his mathematical insight and for suggesting the problem investigated in this paper. Moreover, the author is indebted to the referee for his generous and numerous comments that helped improve the quality of this work.

\section{Background and motivation}
Let $H$ be a Hilbert space and $T$ a bounded linear operator on $H$, which we indicate by $T\in\B(H)$. Assume that $T$ is a completely non-unitary contraction. Let $H^\infty$ be the Hardy space of bounded holomorphic functions on the unit disc $\D$. The Sz.-Nagy--Foias $H^\infty$ \textit{functional calculus} then provides a contractive algebra homomorphim $\Phi: H^\infty\to \B(H)$ such that $\Phi(p)=p(T)$ for every polynomial $p$. Moreover, $\Phi$ is continuous when $H^\infty$ and $\B(H)$ are given their respective weak-star topologies. We will write $\Phi(u)=u(T)$ for $u\in H^\infty$. The contraction $T$ is said to be \textit{of class} $C_0$ if  $\Phi$ has non-trivial kernel. It is known that $\ker \Phi=m_T H^\infty$ for some inner function $m_T$ called the \textit{minimal function} of $T$. The minimal function is uniquely determined up to a scalar of absolute value one.

For any inner function $\theta\in H^\infty$, the space $H(\theta)=H^2\ominus \theta H^2$ is closed and invariant for $S^*$, the adjoint of the shift operator $S$ on $H^2$. The operator $S(\theta)$ defined by $S(\theta)^*=S^*|(H^2\ominus \theta H^2)$ is called a \textit{Jordan block}; it is of class $C_0$ with minimal function $\theta$.

A vector $x\in H$ is said to be \textit{cyclic} for $T\in \B(H)$ if the linear manifold generated by $\{T^n x: n\geq 0\}$ is dense in $H$. If an operator has a cyclic vector, it is called \textit{multiplicity-free}. A bounded linear operator $X:H\to H'$ is called a \textit{quasiaffinity} if it is injective and has dense range. The following is Theorem 3.2.3 of \cite{bercOTA}, its conclusion is summarized by saying that $T$ is \textit{quasisimilar} to $S(m_T)$. 

\begin{theorem}
Let $T\in \B(H)$ be a multiplicity-free operator of class $C_0$. Then, there exist quasiaffinities $X: H\to H(m_T)$
and $Y:H(m_T)\to H$ with the property that $XT=S(m_T)X$ and $TY=YS(m_T).$
\end{theorem}
More details about all of the above background material can be found in \cite{bercOTA}.

Throughout the paper, we use the following notation for Blaschke factors: for $\lambda\in \D,\lambda\neq 0$ we set
$$
b_{\lambda}(z)=\frac{\overline{\lambda}}{|\lambda|}\frac{\lambda-z}{1-\overline{\lambda}z},
$$
and we also set $b_0(z)=z$.

We now give motivation for the main problem that we address. Let $\{\lambda_j\}_j\subset \D$ be a \textit{Carleson sequence,} that is
a sequence satisfying (\ref{unifsep}). Set $b(z)=\prod_{j\in \N}b_{\lambda_j}(z).$
The classical Carleson interpolation theorem (see
\cite{carleson}, \cite{garnett}) implies that the map
$$
H^\infty/ b H^\infty \to \ell^\infty
$$
$$
u+b H^\infty\mapsto \{u(\lambda_j)\}_j
$$
is a bounded algebra isomorphism. On the other hand, a consequence
of the commutant lifting theorem (see \cite{Sarason}, \cite{NF1}, \cite{NF2}) is that the map
$$
H^\infty/b H^\infty\to \{S(b)\}'
$$
$$
u+b H^\infty\mapsto u(S(b))
$$
is an isometric algebra isomorphism. Thus, we see that there exists a bounded algebra isomorphism between $\{S(b)\}'$ and $\ell^\infty.$
Consider now the operator
$$
T=\bigoplus_{j\in \N}S\left(b_{\lambda_j}\right).
$$
It is obvious that $T$ is of class $C_0$ with minimal function $b$. Moreover, any vector $\bigoplus_{j\in \N}h_j\in \bigoplus_{j\in \N}H(b_{\lambda_j})$ with $h_j\neq 0$ for every $j\in \N$ is cyclic for $T$. Therefore, $T$ is quasisimilar to $S(b)$. If we denote by $\id_X$ the identity operator on a space $X$, then we have 
$
S(b_{\lambda_j})=\lambda_j \id_{\C},
$
and thus $$\{T\}'=\left\{\bigoplus_{j\in \N}a_j:\{a_j\}_j\in \ell^\infty\right\}\cong \ell^\infty,$$  so in fact $\{S(b)\}'$ and $\{T\}'$ are boundedly isomorphic as algebras. We thus see that $T$ and $S(b)$ share properties beyond what is guaranteed by mere quasisimilarity. We set out to investigate this phenomenon more carefully.

\section{Blaschke product with multiplicity one}
Let $\{\lambda_j\}_j\subset \D$ satisfy the Carleson condition (\ref{unifsep}), and set 
$$
b(z)=\prod_{j\in \N}b_{\lambda_j}(z).
$$
By the classical Carleson interpolation theorem, for every subset
$A\subset \N$, we can find $\phi_A\in H^\infty$ such that
$$
\phi_A(\lambda_j)=\begin{cases}1 \text{ if } j\in A\\
                                    0 \text{ if } j \notin A
                                    \end{cases}
$$ with $\|\phi_A\|\leq C$,where $C$ is the  so-called constant of interpolation (see \cite{garnett} p.276) and is independent of the set $A$.

\begin{lemma}\label{group}
Let $\Psi:H^\infty \to \B(H)$ be a bounded algebra homomorphism such that $b H^\infty \subset\ker \Psi$. For
every $A\subset \N$, define
$$
g_A=\Psi(\phi_A)-(\id_{H}-\Psi(\phi_A))\in
\B(H).
$$
Then
$$G=\{g_A:A\subset \N\}\subset \B(H)$$ is an abelian group under
multiplication.
\end{lemma}
\begin{proof}
We first need to check that $g_A$ is well-defined since the function $\phi_A$ is not uniquely determined. Assume that $\phi_1$ and $\phi_2$ are two candidates for $\phi_A$. Then, $\phi_1-\phi_2$ vanishes at every $\lambda_j$, whence $\phi_1-\phi_2=bf$ for some $f\in H^\infty$. Consequently, $\Psi(\phi_1-\phi_2)=\Psi(bf)=0$ since $\Psi$ vanishes on $bH^\infty$, and thus $g_A$ is well-defined. A straightforward calculation now yields that $\Psi(\phi_A)$ is idempotent and that $g_A$ is an invertible operator. Moreover, $g_A
g_B=g_{(A\cap B)\cup (A^c\cap B^c)}$.
\end{proof}

\begin{proposition}\label{directsum}
Let $T\in \B(H)$. Assume that there exists a bounded algebra
homomorphism $\Psi: H^\infty\to \B(H)$ such that $\Psi(p)=p(T)$ for
every polynomial $p$, which is also continuous when $H^\infty$ is given the weak-star topology and
$\B(H)$ the weak operator topology . Furthermore,
assume that $H=\bigvee_{j\in \N}H_j$ where $H_j=\ker (T-\lambda_j)$
and $\{\lambda_j\}_j\subset \D$ is a Carleson sequence.
Then there exists an invertible operator $X\in \B(H)$ such that the subspaces $\{XH_j\}_{j\in \N}$ are mutually othogonal and $XTX^{-1}=\bigoplus_{j\in \N}\lambda_j \id_{XH_j}$.
\end{proposition}

\begin{proof}
Since $H=\bigvee_{j\in \N}H_j$, it is easy to check that
$bH^\infty\subset \ker \Psi$. Moreover, the fact that $\Psi$ is
bounded implies $$\|g_A\|\leq 2 \|\Psi\| \|\phi_A\|+1$$ for every
$A\subset \N$, and using the fact that $\|\phi_A\|_{H^\infty}\leq
C$, we see that the inclusion $\iota:G\to \B(H)$ is strongly
continuous and uniformly bounded. By Dixmier's theorem (see
\cite{dixmier} or Theorem 9.3 of \cite{paulsen}), there exists a bounded invertible
operator $X\in \B(H)$ such that $Xg_AX^{-1}\in \B(H)$ is unitary for
every $A\subset \N$.  Now, the equations
$$(Xg_AX^{-1})^*(Xg_AX^{-1})=\id_H=(Xg_AX^{-1})(Xg_AX^{-1})^*$$
along with the definition of $g_A$ easily yield that
$X\Psi(\phi_A)X^{-1}$ must be normal. Hence, $X\Psi(\phi_A)X^{-1}$
is a self-adjoint projection. Notice that for $j\neq k$,  $\phi_{\{j\}}\phi_{\{k\}}=0$ so that
$\sum_{j\in \N}X\Psi(\phi_{\{j\}})X^{-1}$ converges strongly to an element in $\B(H)$.

For each $j\in \N$, choose a sequence $\{p_{j,n}\}_n$ of polynomials
such that $p_{j,n}\to \phi_{\{j\}}$ in the weak-star topology as
$n\to \infty$. By the continuity property of $\Psi$, we have
$$\Psi(p_{j,n})\to \Psi(\phi_{\{j \}})$$
in the weak operator topology as $n\to \infty$. In particular, we see
that for $x\in H_k$ and $y\in H$, we have
$$
\langle \Psi(\phi_{\{j\}})x,y \rangle=\lim_{n\to \infty}\langle
\Psi(p_{j,n})x,y \rangle=\lim_{n\to \infty}\langle (p_{j,n}(T)|H_k)
x,y \rangle$$ $$=\lim_{n\to \infty} p_{j,n}(\lambda_k)\langle
x,y\rangle=\phi_{\{j\}}(\lambda_k)\langle x,y\rangle
$$
whence $\Psi(\phi_{\{j\}})x=\phi_{\{j\}}(\lambda_k)x$ for every $x\in H_k$ and every $k\in \N$. By choice of $\phi_{\{j\}}$, this shows that $\Psi(\phi_{\{j\}})=\text{Id}$ on $H_j$, $\Psi(\phi_{\{j\}})=0$ on $H_k$ for $k\neq j$  and
$$
Tx=\lambda_k x=\sum_{j\in \N}\lambda_j\Psi(\phi_{\{j\}})x
$$
for every $x\in H_k$ and every $k\in \N$. Since $H=\bigvee_{k\in \N}H_k$, we get that
$$
XTX^{-1}=\sum_{j\in \N} \lambda_j X\Psi(\phi_{\{j\}})X^{-1}.
$$
Recall now that every $X\Psi(\phi_{\{j\}})X^{-1}$ is a self-adjoint projection. We conclude that $XTX^{-1}$ is normal, so that $XH_j\perp XH_k$ for $j\neq k$ and $XTX^{-1}=\bigoplus_{j\in\N}\lambda_j \id_{XH_j}$.
\end{proof}

\begin{corollary}\label{similardiag}
Let $\{\lambda_j\}_j\subset \D$ be a Carleson sequence and let $T\in \B(H)$ be an operator of class $C_0$ with minimal function $b(z)=\prod_{j\in \N}b_{\lambda_j}(z)$. Then $T$ is similar to $\bigoplus_{j\in \N}\lambda_j \id_{H_j}$ where $H_j=\ker (T-\lambda_j)$ for every $j\in \N$.
\end{corollary}
\begin{proof}
Notice that $b$ is the least common inner multiple of the family $\{b_{\lambda_j}\}_j$. By Theorem 2.4.6 in \cite{bercOTA}, we can thus write
$$
H=\ker b(T)=\bigvee_{j\in \N}\ker b_{\lambda_j}(T)=\bigvee_{j\in \N}H_j.
$$
Apply now Proposition \ref{directsum} to $T$ with $\Psi$ being the usual Sz.-Nagy--Foias $H^\infty$ functional calculus. We see that there exists an invertible operator $X\in \B(H)$ such that the subspaces $\{XH_j\}_{j\in \N}$ are mutually othogonal and $XTX^{-1}=\bigoplus_{j\in \N}\lambda_j \id_{XH_j}$. 
Let $H'=\bigoplus_{j\in \N}H_j$ be the external orthogonal sum
of the subspaces $H_j\subset H$, and let $\iota_j:H_j\to H'$ denote the
canonical inclusion. Define $X':H'\to H$ as
$X'|\iota_j(H_j)=X|H_j$. Using that $XH_j\perp XH_k$ for $j\neq k$, it is easy to verify that $X'$ is bounded and
invertible, and that it establishes a similarity between $\bigoplus_{j\in
\N}\lambda_j \id_{H_j}$ and $\bigoplus_{j\in \N}\lambda_j \id_{XH_j}.$ The proof is complete.
\end{proof}

\section{Characterization of Carleson sequences}

We now want to give a characterization of the sequences
$\{\lambda_j\}_j\subset \D$ which satisfy the Carleson
condition (\ref{unifsep}). This is done in Theorem \ref{equivalence2}. 
Assume that $\{\lambda_j\}_j\subset \D$ is a
Blaschke sequence of distinct points. Let 
$$
b(z)=\prod_{j\in
\N}b_{\lambda_j}(z).
$$

\begin{lemma}\label{closedrangeD}
Let $D=\bigoplus_{j\in \N}S(b_{\lambda_j})$ and $u\in H^{\infty}$ be an inner divisor of $b$. Then $u(D)$ has closed range if and
only if
$$\inf\{|u(\lambda_j)|:u(\lambda_j)\neq 0 \}>0.$$
\end{lemma}
\begin{proof}
It is easy to verify that $D=\bigoplus_{j\in \N}\lambda_j\text{Id}_{\C}$. Let $F=\{j\in \N: u(\lambda_j)\neq 0\}$. Then $u(D)=\bigoplus_{j\in \N}u(\lambda_j)$ has closed
range if and only if $\bigoplus_{j\in F}u(\lambda_j)$ has closed range. But
$\bigoplus_{j\in F}u(\lambda_j)$
is clearly injective and has dense range, so
that $\bigoplus_{j\in F}u(\lambda_j)$ has closed range if and only $\bigoplus_{j\in F}u(\lambda_j)$ is
invertible, which in turn is equivalent to the condition of the
lemma.
\end{proof}

For $n\in \N$, put
$$
u_n(z)=\frac{b(z)}{b_{\lambda_n}(z)}=\prod_{j\neq
n}\left(\frac{\overline{\lambda_j}}{|\lambda_j|}\frac{\lambda_j-z}{1-\overline{\lambda_j}z}\right).
$$
Moreover, we will need a notation for partial products of an infinite product. Given $f(z)=\prod_{j}f_j(z)$, denote by $f(z;m)$ the $m$-th partial product $\prod_{j=1}^m
f_j(z)$.

\begin{lemma}\label{construction}
Assume that $u_n(\lambda_n)\to 0$ as $n\to \infty$. Then there exists an inner divisor $u$ of $b$ such that
$$\inf\{|u(\lambda_j)|:u(\lambda_j)\neq 0 \}=0.$$

\end{lemma}
\begin{proof}
Define $\alpha_k=1+2^{-(k-1)}$. Choose $n_1\in \N$ with the property that $|u_{n_1}(\lambda_{n_1})|<2^{-2}$. Choose $m_1>n_1$ such that $|u_{n_1}(\lambda_{n_1};m_1)|<2^{-1}$ and note that $u_{n_1}(\lambda_{n_1};m_1)\neq 0$. Set $w_1(z)=u_{n_1}(z;m_1).$ Next, choose $n_2>m_1$ with the property that $|u_{n_2}(\lambda_{n_2})|<2^{-3}$. Note that $|b_{\lambda}(z)|\to 1$ as $|\lambda|\to 1$ for every $z\in \D$, so we can also require that $|b_{\lambda_{n_2}}(\lambda_{n_{1}})|\geq \alpha_2^{-1}$. Define $$v_2(z)=\prod_{j\neq n_1,n_2} \left(\frac{\overline{\lambda_j}}{|\lambda_j|}\frac{\lambda_j-z}{1-\overline{\lambda_j}z}\right).$$ Notice that $v_2 b_{\lambda_{n_1}}=u_{n_2}$ and $v_2 b_{\lambda_{n_2}}=u_{n_1}$. Hence,
$$|v_2(\lambda_{n_2})|=\left|\frac{u_{n_2}(\lambda_{n_2})}{b_{\lambda_{n_1}}(\lambda_{n_2})}\right|\leq \frac{\alpha_2}{2^3}$$ and
$$|v_2(\lambda_{n_1})|=\left|\frac{u_{n_1}(\lambda_{n_1})}{b_{\lambda_{n_2}}(\lambda_{n_1})}\right|\leq \frac{\alpha_2}{2^2}$$
where we used the fact that $|b_x(y)|=|b_y(x)|$. Choose
$m_2>n_2$ such that $$|v_2(\lambda_{n_2};m_2)|\leq
\frac{\alpha_2}{2^2}$$ and $$|v_2(\lambda_{n_1};m_2)|\leq
\frac{\alpha_2}{2^1}.$$ Set $w_2(z)=v_2(z;m_2).$

Assume that we have chosen $n_1,\ldots,n_{k-1},m_1,\ldots,
m_{k-1}\in \N$ with $m_j>n_j>m_{j-1}>n_{j-1}$ for $2\leq j \leq k-1$ and  $v_1,\ldots,
v_{k-1},w_1,\ldots,w_{k-1}$ as above with the property that
$$|v_j(\lambda_{j})|\leq \frac{1}{2^{j+1}}\prod_{\mu=2}^{j}\alpha_{\mu}$$ and
$$|v_j(\lambda_{p})| \leq \frac{1}{2^{p+1}}\prod_{\mu=2}^{j}\alpha_{\mu}$$ for every $1\leq p<j\leq k-1$. We will construct $w_{k}$.  Choose $n_k>m_{k-1}\in \N$ such that $|u_{n_k}(\lambda_{n_k})|<2^{-(k+1)}$ and
$$|b_{\lambda_{n_k}}(\lambda_{n_{j}})|\geq \frac{1}{\alpha_k}$$ for every $j=1,\ldots,k-1$. Define $$v_k(z)=\prod_{j\neq n_1,\ldots,n_k} \left(\frac{\overline{\lambda_j}}{|\lambda_j|}\frac{\lambda_j-z}{1-\overline{\lambda_j}z}\right).$$ Notice that
$v_k b_{\lambda_{n_k}}=v_{{k-1}}$ and $v_k \prod_{j=1}^{k-1}b_{\lambda_{n_j}}=u_{n_k}$. Hence,
$$|v_k(\lambda_{n_k})|=\left|\frac{u_{n_k}(\lambda_{n_k})}{\prod_{\mu=1}^{k-1}b_{\lambda_{n_{\mu}}}(\lambda_{n_k})}\right|\leq \frac{\alpha_k^{k-1}}{2^{k+1}}\leq \frac{1}{2^{k+1}}\prod_{\mu=2}^{k}\alpha_{\mu}$$
and
$$|v_k(\lambda_{n_j})|=\left|\frac{v_{{k-1}}(\lambda_{n_{j}})}{b_{\lambda_{n_k}}(\lambda_{n_{j}})}\right|\leq  \frac{1}{2^{j+1}}\prod_{\mu=2}^{k}\alpha_{\mu}$$
for $j=1,\ldots,k-1$. Choose $m_k\in \N$ such that $$|v_k(\lambda_{n_k};m_k)|\leq \frac{1}{2^{k}}\prod_{\mu=2}^{k}\alpha_{\mu}$$ and $$|v_k(\lambda_{n_j};m_k)|\leq
\frac{1}{2^{j}}\prod_{\mu=2}^{k}\alpha_{\mu}$$ for $j=1,\ldots,k-1$.
Set $w_k(z)=v_k(z;m_k).$

Finally, define $u=\lim_k w_k$. It is straightforward to check that for every $p\in \N$ the function $u$ defines an inner
divisor of $b/b_{\lambda_{n_p}}$, whence $u(\lambda_{n_p})\neq 0$ and $u$ is an inner divisor of $b$.  It remains to show that
$u(\lambda_{n_p})\to 0$ as $p\to \infty$. By definition, we have
$$|u(\lambda_{n_p})|=\lim_k |w_k(\lambda_{n_p})|\leq \frac{1}{2^p}\lim_k
\prod_{\mu=2}^k \alpha_{\mu},$$ so that $|u(\lambda_{n_p})|\to 0$ as $p\to \infty$ since
$\prod_{\mu} \alpha_{\mu} $ is convergent.
\end{proof}

\begin{proposition}\label{equivalence1}
Let $\{\lambda_j\}_j\subset \D$ be a Blaschke sequence of distinct points. Let $D=\bigoplus_{j\in \N}S(b_{\lambda_j})$. Then
$\{\lambda_j\}_j\subset \D$ is a Carelson sequence if and only if
$u(D)$ has closed range for every inner divisor $u\in H^{\infty}$  of
$b$.
\end{proposition}
\begin{proof}
Assume that $\{\lambda_j\}_j$ is not a Carleson sequence. This means that
$$\inf_{k\in\N}\left| u_k(\lambda_k)\right|=\inf_{k\in\N}\prod_{j\neq k}\left| \frac{\lambda_j-\lambda_k}{1-\overline{\lambda_j}\lambda_k}\right|=0$$
so we can find a sequence $\{u_{k_n}\}_n$ such that $u_{k_n}(\lambda_{k_n})\to
0$ as $n\to \infty$. By Lemma \ref{closedrangeD} and
\ref{construction}, we see that there exists some inner divisor $u$
of $b$ for which $u(D)$ doesn't have closed range.

Conversely, assume that $\{\lambda_j\}_j$ is a Carleson sequence, and
let $u\in H^\infty$ be an inner divisor of $b$. Assume that
$u(\lambda_j)\neq 0$. Since $u$ divides $b$, by definition of $u_j$
we see that there must exist an inner function $\psi\in H^\infty$
such that $u_j=\psi u$. In turn, this implies that
$|u(\lambda_j)|\geq |u_j(\lambda_j)|$. By assumption, we know that
$\inf_{j}|u_j(\lambda_j)|>0$, so that
$$\inf\{|u(\lambda_j)|:u(\lambda_j)\neq 0 \}>0.$$ By virtue of Lemma
\ref{closedrangeD}, this completes the proof.
\end{proof}
We can now prove the main result of this section, which gives another interpretation of the classical Carleson interpolation theorem.

\begin{theorem}\label{equivalence2}
Let $\{\lambda_j\}_j\subset \D$ be a Blaschke sequence of distinct points and $b=\prod_{j\in \N}b_{\lambda_j}$ be the corresponding Blaschke product. The following statements are
equivalent:
\begin{enumerate}[(i)]
\item $\{\lambda_j\}_j\subset \D$ is a Carleson sequence 
\item every operator $T$ of class $C_0$ with minimal function $b$ is similar to $\bigoplus_{j\in \N}\lambda_j \id_{H_j}$ where $H_j=\ker (T-\lambda_j)$ for every $j\in \N$
\item every multiplicity-free operator $T$ of class $C_0$ with minimal function $b$ is similar to $S(b)$
\item $u(T)$ has closed range for every multiplicity-free operator $T$ of class $C_0$ with minimal function $b$ and every inner divisor $u\in H^\infty$ of $b$
\item $u(D)$ has closed range for every inner divisor $u\in H^\infty$ of $b$, where  $D=\bigoplus_{j\in \N}S(b_{\lambda_j})$.
\end{enumerate}
\end{theorem}
\begin{proof}
Corollary \ref{similardiag} shows that (i) implies (ii).

Assume that (ii) holds  and let $T$ be a multiplicity-free operator of class $C_0$ with minimal function $b$. Set $H_j=\ker (T-\lambda_j)$ and $K_j=\ker (S(b)-\lambda_j)$ for every $j\in \N$, and choose $X:H\to H(b)$ and $Y:H(b)\to H$ quasiaffinities such that $XT=S(b)X$ and $ TY=YS(b)$. A routine verification shows that $K_j$ has dimension $1$. Now, $XH_j=X\ker (T-\lambda_j)\subset \ker (S(b)-\lambda_j)=K_j$ and $YK_j=Y\ker (S(b)-\lambda_j)\subset \ker (T-\lambda_j)=H_j$. Since $X$ and $Y$ are injective, we find that $H_j$ also has dimension $1$ and $XH_j=K_j$. Corollary \ref{similardiag} then implies that both $T$ and $S(b)$ are similar to $\bigoplus_{j\in \N}\lambda_j \id_{\C}$, which in turn implies (iii).

Assume that (iii) holds. Then $\phi(T)$ is similar to $\phi(S(b))$
for every $\phi\in H^\infty$. Now, it is a  classical fact that
$u(S(b))$ is a partial isometry for every inner
divisor $u\in H^\infty$ of $b$ (see problem 11 p.43 of \cite{bercOTA}), so that
$u(T)$ has closed range for every inner divisor of
$u\in H^\infty$  $b$, which is (iv).

Assume that (iv) holds. As was noted in Section 2, $D$ is a multiplicity-free operator of class $C_0$ with
minimal function $b$. Assertion (v) obviously follows.

Finally, the fact that (v) implies (i) follows from Proposition
\ref{equivalence1}.
\end{proof}

\section{Blaschke product with bounded multiplicity}
We are now ready to address the main question in the case where the
multiplicities of the roots of the minimal function are bounded. As
before, let $\{\lambda_j\}_j\subset \D$ be a Blaschke
sequence of distinct points. Throughout this section we assume that the sequence
$\{\lambda_j\}_j$ satisfies the Carleson condition (\ref{unifsep}). Let $\{m_j\}_j\subset
\N$ be a bounded sequence, set $M=\sup_j m_j$ and define
$$
\theta(z)=\prod_{j=1}^\infty\ b^{m_j}_{\lambda_j}(z).
$$ 
By the interpolation theorem for germs of holomorphic functions (from \cite{Vas1},\cite{Vas2}, see also Chapter 9 Section 4 of \cite{NikShift}), for every $u\in H^\infty$ there exists a
function $\widehat{u}\in H^\infty$ with the property that
$$
\widehat{u}(\lambda_j)=u(\lambda_j)
$$ 
and 
$$
\frac{d^p\widehat{u}}{dz^p}(\lambda_j)=0
$$ 
for every $j\in \N$, $1\leq p\leq
M$. Moreover, we have that $$\|\widehat{u}\|_{H^\infty}\leq
C\|\{u(\lambda_j)\}_j\|_{\ell^\infty}\leq C\|u\|_{H^\infty}$$ for
some constant $C>0$ independent of $u$ (this follows from
the open mapping theorem). Denote by $\chi$ the identity function on $\D$.

\begin{lemma}\label{Delta}
Let $T\in \B(H)$ be an operator of class $C_0$ with minimal
function $\theta$. Then the map
$$
\Psi: H^\infty \to \B(H)
$$
$$
u\mapsto \widehat{u}(T)
$$
is a bounded algebra homomorphism such that
$\Psi(p)=p(\widehat{\chi}(T))$ for every polynomial $p$. Moreover, $\Psi$ is
continuous when $H^\infty$ is given the weak-star topology and
$\B(H)$ the weak operator topology. Finally, 
$\ker(\widehat{\chi}(T)-\lambda_j)
=\ker b_{\lambda_j}^{m_j}(T)$ for every $j\in \N$ and $H=\bigvee_{j\in \N}\ker b_{\lambda_j}^{m_j}(T)$.
\end{lemma}

\begin{proof}
Arguing as in the proof of Corollary \ref{similardiag}, we see that 
$$
H=\bigvee_{j\in\N}\ker b_{\lambda_j}^{m_j}(T).
$$
Set
$H_j=\ker b_{\lambda_j}^{m_j}(T).$ A routine verification shows that $u(T)|H_j=u(T|H_j)$ for
every $j\in \N$ and $u\in H^\infty$. Notice now that
$\widehat{u}-u(\lambda_j)$ has a zero of order $M$ at $\lambda_j$,
so that $b_{\lambda_j}^{m_j}$ divides $\widehat{u}-u(\lambda_j)$, whence $\widehat{u}(T)=u(\lambda_j)$ on $H_j$. In particular, we see that $H_j\subset\ker(\widehat{\chi}(T)-\lambda_j)$. Choose $x\in \ker(\widehat{\chi}(T)-\lambda_j)$. The restriction of $T$ to $\bigvee_{n=0}^\infty T^n x$ is still of class $C_0$ (see \cite{bercOTA}) and hence it has a minimal function, say $m_x$. Then $m_x$ divides the greatest common inner divisor of $\theta$ and $\widehat{\chi}-\lambda_j$, and thus $m_x$ divides $b_{\lambda_j}^{m_j}$, which implies that $x\in H_j$. In other words,  $H_j=\ker b_{\lambda_j}^{m_j}(T)=\ker(\widehat{\chi}(T)-\lambda_j)$.

Let us now show that $\Psi$ is well-defined. Assume that $\widehat{u_1}$ and $\widehat{u_2}$ are two bounded holomorphic functions on $\D$ having the same interpolating property as $\widehat{u}$. Then, $\widehat{u_1}-\widehat{u_2}$ has a zero of order $M$ at every $\lambda_j$, and thus $\theta$ divides $\widehat{u_1}-\widehat{u_2}$. But $\theta$ is the minimal function of $T$, so that $(\widehat{u_1}-\widehat{u_2})(T)=0$ and $\Psi$ is well-defined.

It is easily checked that $\Psi$ is bounded:
$$
\|\Psi(u)\|=\|\widehat{u}(T)\|\leq \|\widehat{u}\|\leq C\|u\|.
$$
Choose now a sequence
$\{u_n\}_n\subset H^\infty$ such that $u_n\to u$ in the weak-star
topology, i.e. $u_n(\lambda)\to u(\lambda)$ for every $\lambda\in
\D$, and $\sup_n \|u_n\|<\infty$. Since $\Psi$ is bounded and $H=\bigvee_{j\in \N}H_j$, to check
the desired continuity property it suffices to verify that
$$
\langle \Psi(u)x,y \rangle=\lim_{n\to \infty}\langle \Psi(u_n)x,y
\rangle
$$ for any $y\in H$ and $x\in H_j, j\in \N$. Let us then pick
$x\in H_j, y\in H$. We already established that $\widehat{v}(T)=v(\lambda_j)$ on
$H_j$ for every $v\in H^\infty$, so we see that
$$
\langle \Psi(u)x,y\rangle=\langle \widehat{u}(T)x,y \rangle=\langle
u(\lambda_j)x,y \rangle=\lim_{n\to \infty}\langle
u_n(\lambda_j)x,y\rangle$$ $$=\lim_n \langle \widehat{u_n}(T)x,y
\rangle=\lim_{n\to \infty}\langle \Psi(u_n)x,y\rangle$$ and
$\Psi$ has the announced continuity property.

Finally, if $p$ is a polynomial and $x\in H_j$, then
$$
\Psi(p)x=\widehat{p}(T)x=p(\lambda_j)x=p(\widehat{\chi}(T)|H_j)x=p(\widehat{\chi}(T))x
$$
and, again, since $H=\bigvee_{j\in \N}H_j$ we
conclude that $\Psi(p)=p(\widehat{\chi}(T))$. The homomorphism
property is verified in a similar fashion.
\end{proof}

\begin{proposition}\label{orthogonal}
Let $T\in\B(H)$ be an operator of class $C_0$ with minimal
function
$$
\theta(z)=\prod_{j=1}^\infty\ b^{m_j}_{\lambda_j}(z)
$$ where $\{\lambda_j\}_j\subset \D$ is a Carleson sequence and $\sup_{j}m_j=M<\infty.$
Then $T$ is similar to  $\bigoplus_{j\in \N}T|H_j$, where $H_j=\ker b_{\lambda_j}^{m_j}(T)$.
\end{proposition}

\begin{proof}
Apply Proposition \ref{directsum} to the operator $\widehat{\chi}(T)$
(see Lemma \ref{Delta}) to get that there exists an invertible operator $R\in
\B(H)$ such that $RH_j\perp RH_k$ for $j\neq k$ where $H_j=\ker
(\widehat{\chi}(T)-\lambda_j)=\ker b_{\lambda_j}^{m_j}(T)$ for $j\in \N$. We may let $A_j=RTR^{-1}|RH_j$ to get that $RTR^{-1}=\bigoplus_{j\in \N}A_j$. Let $H'=\bigoplus_{j\in \N}H_j$ be the external orthogonal sum
of the subspaces $H_j\subset H$, and let $\iota_j:H_j\to H'$ denote the
canonical inclusion. Define $R':H'\to H$ as
$R'|\iota_j(H_j)=R|H_j$. Using that $RH_j\perp RH_k$ for $j\neq k$, it is easy to verify that $R'$ is bounded and
invertible, and that it establishes a similarity between $\bigoplus_{j\in
\N}A_j$ and $\bigoplus_{j\in \N}T|{H_j}.$ The proof is complete.
\end{proof}

If in addition we assume that $T$ is multiplicity-free, then we obtain that $T|H_j$ is similar to $S(b_{\lambda_j}^{m_j})$. Indeed, choose $X:H\to H^2\ominus \theta H^2$ and $Y:H^2\ominus \theta H^2\to H$ quasiaffinities such that $XT=S(\theta)X$ and $TY=YS(\theta)$. Arguing as in the proof of Theorem \ref{equivalence2}, it is easy to see that $XH_j=X\ker b_{\lambda_j}^{m_j}(T)=\ker b^{m_j}_{\lambda_j}(S(\theta))$. In particular, $X|H_j$ establishes a similarity between $T|H_j$ and $S(\theta)|\ker b_{\lambda_j}^{m_j}(S(\theta))$, which in turn is unitarily equivalent to $S(b_{\lambda_j}^{m_j})$ by Proposition 3.1.10 in \cite{bercOTA}. This estalishes the claim. However, this procedure doesn't offer any control over the norm of $(X|H_j)^{-1}$, which turns out to be crucial in the proof of our main result. The following development remedies this matter.

Before proceeding we need to recall some classical results from interpolation theory. The first one is a combination of Lemma 3.2.8, Corollary 3.2.11 and Theorem 3.2.14 of
\cite{nikOP}. It is originally from \cite{Vas2}.

\begin{proposition}\label{Nikequiv}
Let $\{\theta_j\}_{ j\in \N}$ be inner functions such that $\theta=\prod_j
\theta_j$ is their least common inner multiple. The following statements are equivalent:
\begin{enumerate}[(i)]
\item there exists a constant $c_1>0$ independent of $z$ with the
property that $|\theta(z)|\geq c_1 \inf_{j\in \N}|\theta_j(z)|$ for
every $z\in \D$
\item there exists a constant $c_2>0$ with the property that for any
finite set $\sigma \subset \N$ there are functions $f_{\sigma},
g_{\sigma}\in H^\infty$ satisfying
$$
f_{\sigma}\theta_{\sigma}+g_{\sigma}\theta_{\N\setminus \sigma}=1
$$
with $\|f_{\sigma}\|\leq c_2$, where $\theta_{\sigma}=\prod_{j \in
\sigma}\theta_j$
\item $\inf_{\sigma\subset \N}\inf_{z\in
\D}\{|\theta_{\sigma}(z)|+|\theta_{\N\setminus \sigma}(z)|\}>0.$
\end{enumerate}
\end{proposition}
The following is Lemma 3.2.18 from \cite{nikOP}.
\begin{proposition}\label{NikGC}
A sequence $\{\lambda_j\}_j\subset \D$ is a Carleson sequence if and
only if statement (1) above holds with $\theta_j=b_{\lambda_j}$.
\end{proposition}

We can now establish an important technical result.

\begin{lemma}\label{closedrange}
Let $\{\lambda_j\}_j\subset \D $ be a Carleson sequence
and $T=\bigoplus_{j\in \N}T_j$ be an operator of class $C_0$
acting on $H=\bigoplus_{j\in \N}H_j$, with minimal function $\theta$ and
such that $b_{\lambda_j}^{m_j}(T_j)=0$. Assume also that $u(T)$ has
closed range for every inner divisor $u\in H^\infty$ of $\theta$. Then
the operator $\bigoplus_{j\in \N}b_j^{k_j}(T_j)$ has closed range
for any choice of integer sequence $\{k_j\}_j$ satisfying $0\leq
k_j\leq m_j$ for every $j\in\N$.
\end{lemma}
\begin{proof}
By Proposition \ref{NikGC}, we have
$$
\left|\prod_{j\in \N}b_{\lambda_j}(z)\right|\geq c_1 \inf_{j\in
\N}|b_{\lambda_j}(z)|
$$
for some constant $c_1>0$ and every $z\in \D$. But then Proposition \ref{Nikequiv} implies
$$
\inf_{\sigma\subset \N}\inf_{z\in \D}\left\{\left|\prod_{j\in
\sigma}b_{\lambda_j}(z) \right|+\left|\prod_{j\notin
\sigma}b_{\lambda_j}(z) \right|\right\}>0.
$$
Hence,
$$
\inf_{\sigma\subset \N}\inf_{z\in \D}\left\{\left|\prod_{j\in
\sigma}b_{\lambda_j}^M(z) \right|+\left|\prod_{j\notin
\sigma}b^M_{\lambda_j}(z) \right|\right\}>0
$$
and in particular, since $m_j\leq M$ and
$b_{\lambda_j}$ is inner for every $j\in \N$ , we may write
$$
\inf_{\sigma\subset \N}\inf_{z\in \D}\left\{\left|\prod_{j\in
\sigma}b_{\lambda_j}^{m_j}(z) \right|+\left|\prod_{j\notin
\sigma}b^{m_j}_{\lambda_j}(z) \right|\right\}
$$
$$
\geq
\inf_{\sigma\subset \N}\inf_{z\in \D}\left\{\left|\prod_{j\in
\sigma}b_{\lambda_j}^M(z) \right|+\left|\prod_{j\notin
\sigma}b^M_{\lambda_j}(z) \right|\right\}>0.
$$
Applying Proposition \ref{Nikequiv} again, we find that there exists
a constant $c_2>0$ with the property that for every $j\in \N$ there
exists functions $f_j,g_j\in H^\infty$ satisfying
$$f_j \frac{\theta}{b_{\lambda_j}^{m_j}}+g_jb^{m_j}_{\lambda_j}=1$$
and $\|f_j\|\leq c_2$. In particular, we have
$$
\id_{H_j}=f_j(T_j) \frac{\theta}{b_{\lambda_j}^{m_j}}(T_j)+g_j(T_j)b_{\lambda_j}^{m_j}(T_j)=f_j(T_j) \frac{\theta}{b_{\lambda_j}^{m_j}}(T_j)=\frac{\theta}{b_{\lambda_j}^{m_j}}(T_j)f_j(T_j) 
$$ 
or
$$
f_j(T_j)=\left( \frac{\theta}{b_{\lambda_j}^{m_j}}(T_j)\right)^{-1}.
$$ 
Note also that $\|f_j(T_j)\|\leq c_2$.
Now, let $\phi=\prod_{r\in \N}b_{\lambda_r}^{k_r}$. Since $0\leq k_j\leq
m_j$, we can write
$$
\frac{\theta}{b_{\lambda_j}^{m_j}}=\frac{\phi}{b_{\lambda_j}^{k_j}}\psi_j
$$
for some inner function $\psi_j\in H^\infty$ and we get the relation
$$
\frac{\theta}{b_{\lambda_j}^{m_j}}(T_j)=\frac{\phi}{b_{\lambda_j}^{k_j}}(T_j)\psi_j(T_j)=\psi_j(T_j)\frac{\phi}{b_{\lambda_j}^{k_j}}(T_j.
$$
Since $\theta/(b_{\lambda_j}^{m_j})(T_j)$ is invertible, we conclude that
$\phi/(b_{\lambda_j}^{k_j})(T_j)$ is actually invertible. We have
$$
\left(
\frac{\phi}{b_{\lambda_j}^{k_j}}(T_j)\right)^{-1}=\psi_j(T_j)f_j(T_j)
$$
and we infer 
$$
\left\| \left(
\frac{\phi}{b_{\lambda_j}^{k_j}}(T_j)\right)^{-1}\right\|\leq
c_2\|\psi_j(T_j)\|\leq c_2
$$ since $\psi_j$ is inner for every $j\in
\N$. This shows that the operator 
$$
\bigoplus_{j\in
\N}\frac{\phi}{b_{\lambda_j}^{k_j}}(T_j)
$$ 
is invertible, and we may
write
$$
\bigoplus_{j\in \N}b_{\lambda_j}^{k_j}(T_j)=\bigoplus_{j\in
\N}b_{\lambda_j}^{k_j}(T_j)\frac{\phi}{b_{\lambda_j}^{k_j}}(T_j)
\left(\frac{\phi}{b_{\lambda_j}^{k_j}}(T_j)\right)^{-1}
$$
$$
=\left(\bigoplus_{j\in
\N}b_{\lambda_j}^{k_j}(T_j)\frac{\phi}{b_{\lambda_j}^{k_j}}(T_j)
\right)\left( \bigoplus_{j\in
\N}\frac{\phi}{b_{\lambda_j}^{k_j}}(T_j)\right)^{-1}$$
$$=\left(\bigoplus_{j\in \N}\phi(T_j) \right)\left( \bigoplus_{j\in \N}\frac{\phi}{b_{\lambda_j}^{k_j}}(T_j)\right)^{-1}=\phi(T)\left( \bigoplus_{j\in \N}\frac{\phi}{b_{\lambda_j}^{k_j}}(T_j)\right)^{-1}
$$
which completes the proof since $\phi$ is an inner divisor of $\theta$, and consequently $\phi(T)$ has closed range by assumption.
\end{proof}

The last ingredient we need is the following. For $m\in \N$, denote by $\chi^m$ the function $z\mapsto z^m$ on $\D$.

\begin{proposition}\label{nilpotent}
Let $H$ be an $n$-dimensional Hilbert space. Let $N\in \B(H)$ be a
multiplicity-free contraction which is nilpotent of order $n$. Then
there exists an invertible operator $X$ such that $XNX^{-1}=S(\chi^n)$ with
the additional property that $\|X\|\leq 1$ and $\|X^{-1}\|\leq \e^{-2(n-1)}$, where $\e>0$ satisfies $\|Nx\|\geq \e \|x\|$ for
every $x\in (\ker N)^{\perp}$.
\end{proposition}
\begin{proof}
There is no loss of generality in assuming that $H=\C^n$. Let $\xi\in
\C^n$ be a cyclic vector for $N$. We may additionnally assume that
$\|N^{n-1}\xi\|=1$. Define now
$$\zeta_{n-1}=N^{n-1}\xi,$$
$$\zeta_{n-2}=N^{n-2}\xi-a_{n-2} N^{n-1}\xi,$$
$$\zeta_{n-3}=N^{n-3}\xi-a_{n-2} N^{n-2}\xi-a_{n-3}N^{n-1}\xi,$$
$$\cdots$$
$$\zeta_{0}=\xi-a_{n-2}N\xi-\ldots-a_{0}N^{n-1}\xi$$
where $a_0,\ldots,a_{n-1}\in \C$ are chosen such that $\langle
\zeta_{n-1},\zeta_k \rangle=0$ for every $0\leq k \leq n-2$. We see that
$\zeta_{0}$ is a cyclic vector for $N$, and hence we may assume that
$N^{n-1}\xi$ is orthogonal to $N^k\xi$ for every $0\leq k \leq n-2$.

Denote by $e_0,\ldots,e_{n-1}$ the usual orthonormal basis of
$\C^n$. Since $\xi$ is cyclic, the vectors $\xi, N\xi, \ldots, N^{n-1}\xi$ form another (possibly non-orthonormal) basis for $\C^n$. Define $X: \C^n\to \C^n$ as $XN^k\xi=e_k$. It is clear then
that $X$ is an invertible linear operator, and we have that
$XNX^{-1}e_{n-1}=XN^n \xi=0$ and $XNX^{-1}e_{k}=XN^{k+1}\xi=e_{k+1}$
for $0\leq k \leq n-2$. In other words, $XNX^{-1}=S(\chi^n)$.

Let us now examine the norm of $X$ and $X^{-1}$. Let
$c_0,\ldots,c_{n-1}$ be arbitrary complex numbers. Using that
$N^{n-1}\xi \perp N^k \xi$ for every $0\leq k \leq n-2$ and that
$\|N^{n-1}\xi\|=1$, we see that
$$
\left\| \sum_{j=0}^{n-1}c_j N^j \xi \right\|^2=\left\|
\sum_{j=0}^{n-2}c_j N^j \xi \right\|^2+|c_{n-1}|^2 .
$$
Now, $\|N\|\leq 1$ by assumption so that we have
\begin{align*}
\left\| \sum_{j=0}^{n-1}c_j N^j \xi \right\|^2&\geq \left\|
\sum_{j=0}^{n-2}c_j N^{j+1} \xi \right\|^2+|c_{n-1}|^2\\
& =\left\|
\sum_{j=0}^{n-3}c_j N^{j+1} \xi \right\|^2+|c_{n-2}|^2+|c_{n-1}|^2.
\end{align*}
Iterating this process yields
$$
\left\| \sum_{j=0}^{n-1}c_j N^j \xi \right\|^2\geq
\sum_{j=0}^{n-1}|c_j|^2 =\left\| \sum_{j=0}^{n-1}c_j e_j \right\|^2.
$$
In other words, $\|X\|\leq 1$. On the other hand, notice that $\ker N=\C
N^{n-1}\xi$ so that
$$
(\ker N)^{\perp}=\bigvee_{j=0}^{n-2}\C N^{j}\xi.
$$
In particular, we see that
$$
\left\| \sum_{j=0}^{n-2}c_j N^j \xi \right\|^2\leq \frac{1}{\e^2}\left\| \sum_{j=0}^{n-2}c_j N^{j+1} \xi \right\|^2
$$
which implies
\begin{align*}
\left\| \sum_{j=0}^{n-1}c_j N^j \xi \right\|^2&=\left\|
\sum_{j=0}^{n-2}c_j N^j \xi \right\|^2 +|c_{n-1}|^2 \\
&\leq
\frac{1}{\e^2}\left\| \sum_{j=0}^{n-2}c_j N^{j+1} \xi
\right\|^2+|c_{n-1}|^2\\
&=\frac{1}{\e^2}\left\| \sum_{j=0}^{n-3}c_j N^{j+1} \xi \right\|^2
+\frac{1}{\e^2}|c_{n-2}|^2+|c_{n-1}|^2 .
\end{align*}
Iterating this calculation yields
$$
\left\| \sum_{j=0}^{n-1}c_j N^j \xi \right\|^2\leq
\sum_{j=0}^{n-1}\frac{|c_j|^2}{\e^{2(n-1-j)}}\leq
\frac{1}{\e^{2(n-1)}}\sum_{j=0}^{n-1}|c_j|^2=\frac{1}{\e^{2(n-1)}}\left\|
\sum_{j=0}^{n-1}c_j e_j \right\|^2.
$$
Thus, $\|X^{-1}\|\leq \e^{-2(n-1)}.$
\end{proof}

We are now ready to prove the main result of the paper.

\begin{theorem}\label{similar2}
Assume that $\{\lambda_j\}_j\subset \D$ is a Carleson sequence
sequence and $\{m_j\}_j\subset \N$ is a bounded sequence. Let $T\in \B(H)$ be
a multiplicity-free operator of class $C_0$ with minimal function $\theta(z)=\prod_{j\in \N}b^{m_j}_{\lambda_j}(z)$. Assume also that $u(T)$ has closed range for every inner divisor $u\in
H^\infty$ of $\theta$. Then $T$ is similar to $\bigoplus_{j\in
\N}S(b_{\lambda_j}^{m_j})$.
\end{theorem}

\begin{proof}
By Proposition \ref{orthogonal}, we may assume that $T=\bigoplus_{j\in \N}T|H_j$ where $H_j=\ker b_{\lambda_j}^{m_j}(T)$. Suppose for the moment that for each $j\in \N$ there exists  an invertible operator $X_j$ such that $X_j b_{\lambda_j}(T|H_j)X_j^{-1}=S(\chi^{m_j})$ and $$
\sup_{j\in
\N}\{\|X_j\|,\|X_j^{-1}\|\}<\infty.
$$
It is known that $b_{\lambda_j}(S(b_{\lambda_j}^{m_j}))$ is
unitarily equivalent to $S(\chi^{m_j})$ (see
problem 2 p.42 of \cite{bercOTA}), so we may assume that
$X_j b_{\lambda_j}(T|H_j)X_j^{-1}=b_{\lambda_j}(S(b_{\lambda_j}^{m_j}))$.
Now, we have $b_{\lambda_j}\circ b_{\lambda_j}=\chi$ on
$\D$, and thus $X_j (T|H_j) X_j^{-1}=S(b_{\lambda_j}^{m_j})$. The bounded operator $X=\bigoplus_{j\in \N}X_j$ is invertible and provides a similarity between $\bigoplus_{j\in \N}T|H_j$ and $\bigoplus_{j\in
\N}S(b_{\lambda_j}^{m_j})$, and the theorem follows.

Hence we are left with the task
of finding such invertible operators $X_j$ for $j\in \N$. By Lemma \ref{closedrange}, we have that $\bigoplus_{j\in \N}b_{\lambda_j}(T|H_j)$ has
closed range, so that there exists $\e>0$ (independent of $j$) with
the property that 
$$
\|b_{\lambda_j}(T|H_j)x\|\geq \e \|x\|
$$ 
for every $x\in (\ker
b_{\lambda_j}(T|H_j))^\perp$ and every $j\in \N$. Note also that the operator $b_{\lambda_j}(S(b_{\lambda_j}^{m_j}))$ is multiplicity-free since it is unitarily equivalent to the Jordan block $S(\chi^{m_j})$, as noted above. The same is true for $b_{\lambda_j}(T|H_j)$ by the remark following Proposition \ref{orthogonal}. Proposition
\ref{nilpotent} then immediately implies that for every $j\in \N$ we can find an invertible operator
$X_j$ intertwining $b_j(T|H_j)$ and $S(\chi^{m_j})$ with the
additional property that $$\sup_{j\in
\N}\{\|X_j\|,\|X_j^{-1}\|\}<\infty$$ since $\sup_{j\in \N}{m_j}<\infty$.
\end{proof}

Finally, we can answer our original question in the case where the
multiplicities $m_j$ are uniformly bounded.

\begin{corollary}\label{reponse}
Assume that $\{\lambda_j\}_j\subset\D$ is a Carleson sequence
sequence and $\{m_j\}_j\subset \N$ is a bounded sequence. Let $T\in \B(H)$ be
a multiplicity-free operator of class $C_0$ with minimal function
$$\theta(z)=\prod_{j=1}^\infty\
\left(\frac{\ol{\lambda_j}}{\lambda_j}\frac{\lambda_j-z}{1-\ol{\lambda}_jz}\right)^{m_j}.$$
Assume that $u(T)$ has closed range for every inner divisor $u\in
H^\infty$ of $\theta$. Then $T$ is similar to $S(\theta)$.
\end{corollary}
\begin{proof}
Both $T$ and $S(\theta)$ satisfy the conditions of Theorem \ref{similar2}, and thus they are similar.
\end{proof}

\end{document}